\documentclass{article}

\usepackage{amsmath}
\usepackage{amssymb}

\newtheorem{thm}{Theorem}

\newtheorem{lem}{Lemma}
\newtheorem{cor}{Corollary}

\title{Size of product of a number and its multiplicative inverse, Moments of L-functions and Exponential Sums}
\author{Tsz Ho Chan}

\begin{document}
\maketitle

\begin{abstract}
In this paper, we study the average size of the product of a number and its multiplicative inverse modulo a prime $p$. This turns out to be related to moments of L-functions and leads to a curious asymptotic formula for a certain triple exponential sum.
\end{abstract}

\section{Introduction and main results}
Let $p$ be a prime number. For any $(a,p) = 1$, let $\overline{a}$ be the positive integer less
than $p$ such that $a \overline{a} \equiv 1 \pmod p$. Of course $a \overline{a}$ can be as small
as $1$ for $a=1$ and as big as $(p-1)^2$ for $a = p-1$. So one can ask on average how big
$a \overline{a}$ is. This leads us to study
\begin{equation} \label{size1}
S := \sum_{a = 1}^{p-1} a \overline{a} = \mathop{\sum_{a = 1}^{p-1}
\sum_{b = 1}^{p-1}}_{a b \equiv 1 \pmod p} a b.
\end{equation}
More generally, one defines
\begin{equation} \label{size2}
S(d) := \mathop{\sum_{a = 1}^{p-1} \sum_{b = 1}^{p-1}}_{a b \equiv d \pmod p}
a b.
\end{equation}
We have
\begin{thm} \label{thm1}
For $(d, p) = 1$,
\[
\mathop{\sum_{a = 1}^{p-1} \sum_{b = 1}^{p-1}}_{a b \equiv d \pmod
p} a b = \frac{p^3}{4} + O(p^{5/2} \log^2 p).
\]
\end{thm}
For a Dirichlet character $\chi$, let $L(s, \chi) = \sum_{n =
1}^{\infty} \frac{\chi(n)}{n^s}$ be the corresponding Dirichlet
$L$-function which has meromorphic continuation over the entire
complex plane. And as a by-product of the proof of Theorem
\ref{thm1}, we have
\begin{cor} \label{cor1}
For $(d,p) = 1$,
\[
\sum_{\chi \neq \chi_0} \overline{\chi}(d) L(0,\chi)^2 \ll p^{3/2}
\log^2 p.
\]
\end{cor}
One can ask how good the error term in Theorem \ref{thm1} is. To
this we consider the mean square error and have
\begin{thm} \label{thm2}
For prime $p$,
\[
\sum_{d = 1}^{p-1} \Bigl| S(d) - \frac{p^2 (p-1)}{4} \Big|^2 =
\frac{5}{144} \frac{p^2 (p^2 - 1)^3}{(p^2 + 1)} + O(p^5 e^{3 \log p
/ \log \log p}).
\]
\end{thm}
This tells us that for some $1 \le d \le p-1$, we have
\[
\Big| S(d) - \frac{p^2 (p-1)}{4} \Big| \gg p^{5/2}.
\]
So the error term in Theorem \ref{thm1} is sharp apart from the
logarithmic factor.

One can consider higher dimensional analogue of (\ref{size2}) by defining
\[
S_k(d) := \mathop{\sum_{a_1 = 1}^{p-1} \sum_{a_2 = 1}^{p-1} ... \sum_{a_k = 1}^{p-1} }_{a_1 a_2
... a_k \equiv d \pmod p} a_1 a_2 ... a_k
\]
and one can prove
\begin{thm} \label{thm4}
For $k \ge 3$ and $(d, p) = 1$,
\[
S_k(d) = \frac{p^{k} (p-1)^{k-1}}{2^k} + O_k(p^{3k/2} (\log p)^k).
\]
\end{thm}
When $k = 3$, one can do slightly better by exponential sum method and get
\begin{thm} \label{thm5}
For $(d, p) = 1$,
\[
S_3(d) = \frac{p^5}{8} + O_k(p^{9/2} (\log p)^2).
\]
\end{thm}
This improvement on the error term may not be very worth doing. But as a by-product of its
proof, we have an interesting result on a triple exponential sum, namely
\begin{thm} \label{thm6}
For $(l, p) =1$,
\[
\sum_{a = 1}^{p-1} \sum_{b = 1}^{p-1} \sum_{c = 1}^{p-1} a b c \; e \Bigl( \frac{l a b c}{p} \Bigr)
= - \frac{p^5}{8} + O(p^{9/2} \log^3 p).
\]
\end{thm}
We will leave the interested readers to derive similar results for exponential sums with more variables.

\bigskip

{\bf Some Notations} Throughout the paper, the symbol $\overline{a}$ stands for the
multiplicative inverse of $a \pmod q$ (i.e. $a \overline{a} \equiv 1
\pmod q$). The notations $f(x) = O(g(x))$, $f(x) \ll g(x)$ and $g(x)
\gg f(x)$ are all equivalent to $|f(x)| \leq C g(x)$ for some
constant $C > 0$. Finally $f(x) = O_\lambda(g(x))$, $f(x)
\ll_\lambda g(x)$ or $g(x) \gg_\lambda f(x)$ mean that the implicit
constant $C$ may depend on $\lambda$.

\section{Some Lemmas}
\begin{lem} \label{lem1}
For $z \neq 1$ and $z^p = 1$, $\sum_{b = 1}^{p-1} b z^b =
\frac{-p}{1 - z}$.
\end{lem}

Proof: As $1 + z + z^2 + ... + z^{p-1} = 0$, one can check directly
that
\[
(\sum_{b = 1}^{p-1} b z^b) (1 - z) = z + z^2 + ... + z^{p-1} - (p-1)
z^p = -1 - (p-1) = -p
\]
which gives the lemma after dividing by $1 - z$.

\begin{lem} \label{lem2}
For $z \neq 1$ and $z^p = 1$, $\sum_{b = 1}^{p-1} \frac{1}{1 - z^b}
= \frac{p-1}{2}$.
\end{lem}

Proof: Notice that $\frac{1}{1 - z} + \frac{1}{1 - \overline{z}} =
\frac{1 - \overline{z} + 1 - z}{(1 - z)(1 - \overline{z})} = \frac{1
- \overline{z} - z + z \overline{z}} {(1 - z)(1 - \overline{z})} =
1$ as $|z| = 1$. Therefore
\[
\sum_{b = 1}^{p-1} \frac{1}{1 - z^b} = \frac{1}{2} \sum_{b =
1}^{p-1} \Bigl( \frac{1}{1 - z^b} + \frac{1}{1 - z^{p-b}} \Bigr) =
\frac{1}{2} \sum_{b = 1}^{p-1} 1 = \frac{p-1}{2}.
\]

\begin{lem} \label{lem3}
For $z \neq 1$, $z^p = 1$ and $1 \le d < p$, $\sum_{b = 1}^{p-1}
\frac{z^{-d b}}{1 - z^b} = \frac{p-1}{2} - d$.
\end{lem}

Proof: Consider
\[
\sum_{b = 1}^{p-1} \frac{1 - z^{-d b}}{1 - z^b} = \sum_{b = 1}^{p-1}
\frac{- z^{-d b}(1 - z^{d b})}{1 - z^b} = - \sum_{b = 1}^{p-1} z^{-d
b} \sum_{j = 0}^{d-1} z^{j b}
\]
\[
= - \sum_{j = 0}^{d-1} \sum_{b = 1}^{p-1} z^{(j - d) b} = - \sum_{j
= 0}^{d-1} (-1) = d.
\]
Therefore by Lemma \ref{lem2},
\[
d = \sum_{b = 1}^{p-1} \frac{1}{1 - z^b} - \sum_{b = 1}^{p-1}
\frac{z^{-d b}}{1 - z^b} = \frac{p-1}{2} - \sum_{b = 1}^{p-1}
\frac{z^{-d b}}{1 - z^b}
\]
which gives the lemma after rearranging terms.

\begin{lem} \label{lem4}
For prime $p$ and $(k,p) = 1$,
\[
\sum_{a = 1}^{p-1} a e \Bigl(\frac{k \overline{a}}{p} \Bigr) \ll
p^{3/2} \log p.
\]
\end{lem}

Proof: By Weil bound on incomplete Kloosterman sum, we have
\[
F(u) := \sum_{a = 1}^{u} e \Bigl(\frac{k \overline{a}}{p} \Bigr) \ll
p^{1/2} \log p
\]
for $1 \le u < p$. Using this and partial summation,
\[
\sum_{a = 1}^{p-1} a e \Bigl(\frac{k \overline{a}}{p} \Bigr) =
\int_{1^-}^{p-1} u dF(u) = (p-1)F(p-1) - \int_{1^-}^{p-1} F(u) du
\ll p^{3/2} \log p.
\]

\begin{lem} \label{lem5}
For $p > 1$,
\[
\sum_{k = 1}^{p-1} \frac{1}{|1 - e(-k/p)|} \ll p \log p.
\]
\end{lem}

Proof: Observe that $|1 - e(-k/p)| \ge |\text{Im} (1 - e(-k/p))| =
|\sin{2 k \pi / p}|$. For $0 \le k < p/4$,  $|\sin{2 k \pi / p}| \ge
k / p$ by observing that the sine function is above the line $y = 2
x / \pi$ for $0 \le x \le \pi / 2$. So
\[
\sum_{k < p/4} \frac{1}{|1 - e(-k/p)|} \le \sum_{k < p/4}
\frac{1}{k/p} \ll p \log p.
\]
Using $\sin(\pi - x) = \sin x$, we have
\[
\sum_{p/4 < k \le p/2} \frac{1}{|1 - e(-k/p)|} \ll p \log p.
\]
Hence
\begin{equation} \label{h1}
\sum_{k = 1}^{p/2} \frac{1}{|1 - e(-k/p)|} \ll p \log p + 1 \ll p
\log p
\end{equation}
where the $1$ may come from the term when $k = p/4$. By complex
conjugation,
\[
\frac{1}{|1 - e(-k/p)|} = \frac{1}{|1 - e(-(p-k)/p)|}.
\]
So from (\ref{h1}),
\begin{equation} \label{h2}
\sum_{k = p/2}^{p-1} \frac{1}{|1 - e(-k/p)|} \ll p \log p
\end{equation}
and the lemma follows from (\ref{h1}) and (\ref{h2}).

\section{Proof of Theorems \ref{thm1} and \ref{thm4} and Corollary \ref{cor1}}
Proof of Theorem \ref{thm1}: We use exponential sum to study
(\ref{size2}). By orthogonality of additive characters,
\[
S(d) = \frac{1}{p} \sum_{a = 1}^{p-1} \sum_{b = 1}^{p-1} a b \sum_{k
= 1}^{p} e \Bigl(\frac{k(d \overline{a} - b)} {p} \Bigr) =
\frac{p(p-1)^2}{4} + \frac{1}{p} \sum_{k = 1}^{p-1} \sum_{a =
1}^{p-1} a e \Bigl(\frac{k d \overline{a}} {p} \Bigr) \sum_{b =
1}^{p-1} b e\Bigl(\frac{-k b} {p} \Bigr)
\]
where $e(u) = e^{2\pi i u}$. Hence, by Lemma \ref{lem1},
\begin{equation} \label{exptemp}
S(d) = \frac{p(p-1)^2}{4} - \sum_{k = 1}^{p-1} \frac{1}{1 - e(-k/p)}
\sum_{a = 1}^{p-1} a e \Bigl(\frac{k d \overline{a}}{p} \Bigr).
\end{equation}
By Lemmas \ref{lem4} and \ref{lem5},
\begin{equation} \label{expresult}
S(d) = \frac{p(p-1)^2}{4} + O \Bigl( p^{3/2} \log p \sum_{k =
1}^{p-1} \frac{1}{|1 - e(-k/p)|} \Bigr) = \frac{p^3}{4} + O(p^{5/2}
\log^2 p).
\end{equation}

Proof of Corollary \ref{cor1}: Another way to study (\ref{size2}) is
through character sums. By orthogonality of Dirichlet characters, we
have
\[
S(d) = \frac{1}{\phi(p)} \sum_{\chi \pmod p} \overline{\chi}(d)
\sum_{a = 1}^{p-1} \sum_{b = 1}^{p-1} a b \chi(a) \chi(b)
\]
\[
= \frac{1}{p-1} \sum_{a = 1}^{p-1} \sum_{b = 1}^{p-1} a b +
\frac{1}{\phi(p)} \sum_{\chi \neq \chi_0} \overline{\chi}(d) \sum_{a
= 1}^{p-1} \sum_{b = 1}^{p-1} a b \chi(a) \chi(b)
\]
\[
= \frac{p^2 (p-1)}{4} + \frac{1}{p-1} \sum_{\chi \neq \chi_0}
\overline{\chi}(d) \Bigl( \sum_{a = 1}^{p-1} a \chi(a) \Bigr)^2.
\]
As $\sum_{a \pmod p} a \chi(a) = -pL(0,\chi)$ (see [\ref{MV}, page 310] and combine with the functional equation for Dirichlet $L$-functions), we have
\begin{equation} \label{charresult}
S(d) = \frac{p^2 (p-1)}{4} + \frac{p^2}{p-1} \sum_{\chi \neq \chi_0}
 \overline{\chi}(d) L(0,\chi)^2.
\end{equation}
Comparing (\ref{charresult}) and (\ref{expresult}), we have
\[
\sum_{\chi \neq \chi_0} \overline{\chi}(d) L(0,\chi)^2 \ll p^{3/2} \log^2 p.
\]

Proof of Theorem \ref{thm4}: The character sum method can be used to study higher dimension
analogue of Theorem \ref{thm1}. By orthogonality of Dirichlet characters, we
have
\[
S_k(d) = \frac{1}{\phi(p)} \sum_{\chi \pmod p} \overline{\chi}(d)
\sum_{a_1 = 1}^{p-1} \sum_{a_2 = 1}^{p-1} ... \sum_{a_k = 1}^{p-1} a_1 a_2 ...
 a_k \chi(a_1) \chi(a_2) ... \chi(a_k)
\]
\[
= \frac{1}{p-1} \sum_{a_1 = 1}^{p-1} \sum_{a_2 = 1}^{p-1} ...
\sum_{a_k = 1}^{p-1} a_1 a_2 ... a_k +
\frac{1}{\phi(p)} \sum_{\chi \neq \chi_0} \overline{\chi}(d) \sum_{a_1 = 1}^{p-1}
\sum_{a_2 = 1}^{p-1} ... \sum_{a_k = 1}^{p-1} a_1 a_2 ... a_k \chi(a_1)
\chi(a_2) ... \chi(a_k)
\]
\[
= \frac{p^{k} (p-1)^{k-1}}{2^k} + \frac{1}{p-1} \sum_{\chi \neq \chi_0}
\overline{\chi}(d) \Bigl( \sum_{a = 1}^{p-1} a \chi(a) \Bigr)^k.
\]
As $\sum_{a \pmod p} a \chi(a) \ll p^{3/2} \log p$ by Polya-Vinogradov inequality
and partial summation, we have
\[
S_k(d) = \frac{p^{k} (p-1)^{k-1}}{2^k} + O_k(p^{3k/2} (\log p)^k)
\]
which gives Theorem \ref{thm4}.

\section{Proof of Theorem \ref{thm2}}
Define
\[
M := \sum_{d = 1}^{p-1} \Bigl| S(d) - \frac{p^2 (p-1)}{4} \Big|^2.
\]
By (\ref{charresult}),
\[
M = \sum_{d = 1}^{p-1} \Bigl| \frac{p^2}{p-1} \sum_{\chi \neq
\chi_0} \overline{\chi}(d) L(0,\chi)^2 \Big|^2
\]
\[
= \frac{p^4}{(p-1)^2} \sum_{\chi_1 \neq \chi_0} \sum_{\chi_2 \neq
\chi_0} L(0,\chi_1)^2 \overline{L(0,\chi_2)^2} \sum_{d = 1}^{p-1}
\overline{\chi_1}(d) \chi_2(d)
\]
\[
= \frac{p^4}{(p-1)} \sum_{\chi_1 \neq \chi_0} |L(0,\chi_1)|^4 =
\frac{p^4}{(p-1)} \mathop{\sum_{\chi_1 \pmod p}}_ {\chi_1(-1) = -1}
|L(0,\chi_1)|^4
\]
by orthogonality of Dirichlet characters and $L(0,\chi) = 0$ when
$\chi(-1) = 1$. Now by $L(0,\chi) = \frac{\tau(\chi)}{\pi} L(1,
\chi)$ and $|\tau(\chi)| = p^{1/2}$,
\[
M = \frac{p^6}{\pi^4 (p-1)} \mathop{\sum_{\chi_1 \pmod p}}_
{\chi_1(-1) = -1} |L(1,\chi_1)|^4 = \frac{5}{144} \frac{p^2 (p^2 -
1)^3}{(p^2 + 1)} + O(p^5 e^{3 \log p / \log \log p})
\]
by Lemma 2 of Zhang [\ref{Z}]. This tells us that for some $1 \le d \le
p-1$, we have
\[
\Big| S(d) - \frac{p^2 (p-1)}{4} \Big| \gg p^{5/2}.
\]
So the error term in (\ref{expresult}) is sharp apart from the
logarithmic factor.

\section{Double exponential sum}
In this section, we want to study
\begin{equation} \label{doubleexp}
D := \sum_{a = 1}^{p-1} \sum_{b = 1}^{p-1} a b \; e \Bigl( \frac{l a
b}{p} \Bigr).
\end{equation}
First, observe that
\[
D = \sum_{d = 1}^{p-1} e \Bigl( \frac{l d}{p} \Bigr)
\mathop{\sum_{a = 1}^{p-1} \sum_{b = 1}^{p-1}}_{a b \equiv d \pmod
p} a b.
\]
By (\ref{exptemp}),
\[
D = - \frac{p(p-1)^2}{4} - \sum_{d = 1}^{p-1} e \Bigl( \frac{l d}{p}
\Bigr) \sum_{k = 1}^{p-1} \frac{1}{1 - e(-k/p)} \sum_{a = 1}^{p-1} a
e \Bigl(\frac{k d\overline{a}}{p} \Bigr).
\]
Now the sums above can be rewritten as
\[
\sum_{k = 1}^{p-1} \frac{1}{1 - e(-k/p)} \sum_{a = 1}^{p-1} a
\sum_{d = 1}^{p-1} e \Bigl(\frac{d (l + k \overline{a})}{p} \Bigr)
\]
\[
= - \sum_{k = 1}^{p-1} \frac{1}{1 - e(-k/p)} \sum_{a = 1}^{p-1} a +
\sum_{k = 1}^{p-1} \frac{1}{1 - e(-k/p)} \sum_{a = 1}^{p-1} a
\sum_{d = 1}^{p} e \Bigl(\frac{d (l + k \overline{a})}{p} \Bigr)
\]
\[ -
\frac{p (p-1)^2}{4} + p \sum_{a = 1}^{p-1} \frac{a}{1 - e(a l / p)}
\]
by Lemma \ref{lem2}. Therefore
\begin{equation} \label{doubleD}
D = - p \sum_{a = 1}^{p-1} \frac{a}{1 - e(a l / p)}.
\end{equation}
\section{Triple exponential sum: Proof of Theorem \ref{thm6}}

In this section, we study
\[
T := \sum_{a = 1}^{p-1} \sum_{b = 1}^{p-1} \sum_{c = 1}^{p-1} a b c \; e \Bigl( \frac{l a b c}{p} \Bigr)
\]
where $0 < l < p$. One can rearrange it as
\[
T = \sum_{c = 1}^{p-1} c \sum_{a = 1}^{p-1} \sum_{b = 1}^{p-1} a b \; e \Bigl( \frac{l c a b}{p} \Bigr)
= -p \sum_{c = 1}^{p-1} c \sum_{a = 1}^{p-1} \frac{a}{1 - e(a c l / p)}
\]
by (\ref{doubleD}). Grouping the sums according to $a c \equiv d \pmod p$, we have
\[
T = - p \sum_{d = 1}^{p-1} \frac{1}{1 - e(d l / p)} \mathop{\sum_{a = 1}^{p-1}
\sum_{c = 1}^{p-1}}_{a c \equiv d \pmod p} a c.
\]
By (\ref{exptemp}),
\[
T = - p \sum_{d = 1}^{p-1} \frac{1}{1 - e(d l / p)} \Bigl[ \frac{p(p-1)^2}{4} -
\sum_{k = 1}^{p-1} \frac{1}{1 - e(-k/p)} \sum_{a = 1}^{p-1} a
e \Bigl(\frac{k d\overline{a}}{p} \Bigr) \Bigr]
\]
\begin{equation} \label{tripleT}
= - \frac{p^2 (p-1)^3}{8} + p \sum_{d = 1}^{p-1} \frac{1}{1 - e(d l / p)} \sum_{k = 1}^{p-1}
\frac{1}{1 - e(-k/p)} \sum_{a = 1}^{p-1} a e \Bigl(\frac{k d\overline{a}}{p} \Bigr)
\end{equation}
by Lemma \ref{lem2}. Theorem \ref{thm6} follows by observing that the above
has a main term $-p^5 / 8$ and an error term $O(p^{9/2} \log^3 p)$ by Lemmas \ref{lem4} and \ref{lem5}.
\section{Proof of Theorem \ref{thm5}}

Now we are ready to study
\[
S_3(d) = \mathop{\sum_{a = 1}^{p-1} \sum_{b = 1}^{p-1} \sum_{c = 1}^{p-1}}_{a b c \equiv d \pmod p} a b c.
\]
By orthogonality of additive characters,
\[
S_3(d) = \frac{1}{p} \sum_{a = 1}^{p-1} \sum_{b = 1}^{p-1} \sum_{c = 1}^{p-1} a b c \sum_{l = 1}^{p} e \Bigl(
\frac{l (a b c - d)}{p} \Bigr) = \frac{p^2 (p-1)^3}{8} + \frac{1}{p} \sum_{l = 1}^{p-1} e \Bigl(-\frac{d l}{p} \Bigr)
\sum_{a = 1}^{p-1} \sum_{b = 1}^{p-1} \sum_{c = 1}^{p-1} a b c e \Bigl( \frac{l a b c}{p} \Bigr).
\]
By (\ref{tripleT}),
\[
S_3(d) = \frac{p^2 (p-1)^3}{8} + \sum_{l = 1}^{p-1} e \Bigl(-\frac{d l}{p} \Bigr) \Bigl[ -
\frac{p (p-1)^3}{8} + \sum_{t = 1}^{p-1} \frac{1}{1 - e(t l / p)} \sum_{k = 1}^{p-1}
\frac{1}{1 - e(-k/p)} \sum_{a = 1}^{p-1} a e \Bigl(\frac{k t\overline{a}}{p} \Bigr) \Bigr]
\]
\[
= \frac{p (p+1) (p-1)^3}{8} + \sum_{t = 1}^{p-1} \sum_{l = 1}^{p-1} \frac{e(-l t
\overline{t \overline{d}}/p)}{1 - e(t l / p)}
\sum_{k = 1}^{p-1} \frac{1}{1 - e(-k/p)} \sum_{a = 1}^{p-1} a e \Bigl(\frac{k t\overline{a}}{p} \Bigr)
\]
\[
= \frac{p (p+1) (p-1)^3}{8} + \sum_{t = 1}^{p-1} \Bigl(\frac{p-1}{2} - \overline{t \overline{d}} \Bigr) \sum_{k = 1}^{p-1}
\frac{1}{1 - e(-k/p)} \sum_{a = 1}^{p-1} a e \Bigl(\frac{k t \overline{a}}{p} \Bigr) =: S_1 + S_2
\]
by Lemma \ref{lem3}. Now
\[
S_2 = \frac{p-1}{2} \sum_{k = 1}^{p-1}
\frac{1}{1 - e(-k/p)} \sum_{a = 1}^{p-1} a \sum_{t = 1}^{p-1} e \Bigl(\frac{k t\overline{a}}{p} \Bigr) +
\sum_{k = 1}^{p-1} \frac{1}{1 - e(-k/p)} \sum_{t = 1}^{p-1} \overline{t \overline{d}} \sum_{a = 1}^{p-1} a
e \Bigl(\frac{k t\overline{a}}{p} \Bigr)
\]
\[
= - \frac{p-1}{2} \sum_{k = 1}^{p-1}
\frac{1}{1 - e(-k/p)} \sum_{a = 1}^{p-1} a + \sum_{k = 1}^{p-1} \frac{1}{1 - e(-k/p)}
\sum_{t = 1}^{p-1} t \sum_{a = 1}^{p-1} a e \Bigl(\frac{k d \overline{t} \overline{a}}{p} \Bigr)
\]
\[
= - \frac{p(p-1)^3}{8} + \sum_{k = 1}^{p-1} \frac{1}{1 - e(-k/p)}
\sum_{t = 1}^{p-1} \sum_{a = 1}^{p-1} a t e \Bigl(\frac{k d \overline{t} \overline{a}}{p} \Bigr)
\]
\[
= - \frac{p(p-1)^3}{8} + \sum_{k = 1}^{p-1} \frac{1}{1 - e(-k/p)}
\sum_{c = 1}^{p-1} e \Bigl(\frac{k d \overline{c}}{p} \Bigr) \mathop{\sum_{t = 1}^{p-1}
\sum_{a = 1}^{p-1}}_{a t \equiv c \pmod p} a t.
\]
By (\ref{exptemp}),
\[
S_2 = - \frac{p(p-1)^3}{8} + \sum_{k = 1}^{p-1} \frac{1}{1 - e(-k/p)}
\sum_{c = 1}^{p-1} e \Bigl(\frac{k d \overline{c}}{p} \Bigr) \Bigl[ \frac{p(p-1)^2}{4} -
\sum_{l = 1}^{p-1} \frac{1}{1 - e(-l/p)} \sum_{a = 1}^{p-1} a
e \Bigl(\frac{l c \overline{a}}{p} \Bigr) \Bigr]
\]
\[
= - \frac{p(p-1)^3}{4} - \sum_{k = 1}^{p-1} \frac{1}{1 - e(-k/p)}
\sum_{c = 1}^{p-1} e \Bigl(\frac{k d \overline{c}}{p} \Bigr)
\sum_{l = 1}^{p-1} \frac{1}{1 - e(-l/p)} \sum_{a = 1}^{p-1} a
e \Bigl(\frac{l c \overline{a}}{p} \Bigr)
\]
\[
= - \frac{p(p-1)^3}{4} - \sum_{k = 1}^{p-1} \frac{1}{1 - e(-k/p)}
\sum_{l = 1}^{p-1} \frac{1}{1 - e(-l/p)} \sum_{a = 1}^{p-1} a \sum_{c = 1}^{p-1}
e \Bigl(\frac{l c \overline{a} + k d \overline{c}}{p} \Bigr)
\]
\[
= - \frac{p(p-1)^3}{4} - \sum_{k = 1}^{p-1} \frac{1}{1 - e(-k/p)}
\sum_{l = 1}^{p-1} \frac{1}{1 - e(-l/p)} \sum_{a = 1}^{p-1} a S(l\overline{a}, k d; p)
\]
where $S(a,b;p)$ is the Kloosterman sum. Using Weil's bound on Kloosterman sum and Lemma \ref{lem5}, we have
\[
S_2 = - \frac{p(p-1)^3}{4} + O(p^{9/2} \log^2 p).
\]
Consequently,
\[
S_3(d) = \frac{p (p-1)^4}{8} + O(p^{9/2} \log^2 p)
\]
which gives Theorem \ref{thm5}.

White Station High School \\
514 S. Perkins Road, \\
Memphis, TN 38117 \\
U.S.A.

\end{document}